\documentclass{article}

\usepackage[utf8]{inputenc}
\usepackage[english]{babel}

\usepackage{graphicx}
\usepackage{amssymb}
\usepackage{amsmath}
\usepackage{tikz-cd}
\usepackage{enumerate}
\usepackage{verbatim}

\usepackage{hyperref}

\newtheorem{theorem}{Theorem}[section]
\newtheorem{corollary}[theorem]{Corollary}

\newtheorem{proposition}[theorem]{Proposition}

\newcommand{\C}{\mathcal{C}}

\begin{document}

\title{Virtual Specht stability for $FI$-modules in positive characteristic}
\author{{Nate Harman} \\
 \textit{\small Department of Mathematics} \\ 
\textit{\small Massachusetts Institute of Technology}  \\
\textit{\small Cambridge, MA, 02139, USA}\\
\texttt{\small nharman@math.mit.edu}}


\maketitle

\vspace{-.9cm}
\begin{abstract}
We define a notion of virtual Specht stability which is a relaxation of the Church-Farb notion of  representation stability for sequences of symmetric group representations.  Using a structural result of Nagpal, we show that $FI$-modules over fields of positive characteristic exhibit virtual Specht stability.
\end{abstract}



\begin{section}{Introduction}

Church and Farb defined the notion of representation stability for a sequence of symmetric group representations \cite{CF} over a field of characteristic zero. Most notably, it has been shown that the sequence of representations defined by a finitely generated $FI$-module over a field of characteristic zero exhibits representation stability \cite{CEF}.

$FI$-modules can be defined over fields of positive characteristic and have been studied in that setting fairly extensively in \cite{CEFN} and \cite{Nagpal}, however so far there has not been a replacement for the notion of representation stability which holds in this context. The purpose of this note is to define a notion of virtual Specht stability which is a relaxation of representation stability that holds for finitely generated $FI$-modules over fields of arbitrary characteristic.

We'll briefly note that theory of finitely generated $FI$-modules has been applied in numerous situations in geometry, topology, and algebra. Theorem \ref{virtstab} will immediately imply that virtual Specht stability holds in many of those applications. In particular virtual Specht stability will hold for the mod-$p$ cohomology of configuration spaces $\text{Conf}_n(M)$, as well as for the homology of congruence subgroups $\Gamma_n(p)$ with coefficients in $\mathbb{F}_p$.  See \cite{CEF} and \cite{CEFN} for a discussion of these and other examples.

\begin{subsection}{$FI$-modules and representation stability}

Let $FI$ denote the category where the objects are finite sets, and morphisms are injections.  An $FI$-module over a commutative ring $k$ is a covariant functor $V$ from $FI$ to the category of modules over $k$. 

Since the set of $FI$-endomorphisms of the set $[n] = \{1,2,\dots,n\}$ is the symmetric group $S_n$, an $FI$-module $V$ can be thought of as a sequence of representations $V_n$ (the image of $[n]$ under the functor) of $S_n$ for each $n$ along with with compatibility maps from $V_n$ to $V_m$ for every injection from $[n]$ to $[m]$. An $FI$-module is said to be finitely generated in degree at most $d$ if all the $V_n$ are finitely generated as $k$-modules and $V_n$ is spanned by the images of $V_d$ under all injections from $[d]$ to $[n]$ for all $n >d$.

$FI$-modules were introduced by Church, Ellenberg, and Farb in \cite{CEF}, where it was shown that over a field of characteristic zero the sequence of symmetric group representations defined by a finitely generated $FI$-module exhibits the phenomenon known as representation stability as defined by Church and Farb in \cite{CF}.

If $\lambda = (\lambda_1, \lambda_2, \dots, \lambda_\ell)$ is a partition, then for $n \ge |\lambda|+\lambda_1$ let $\lambda[n] = (n - |\lambda|, \lambda_1, \lambda_2, \dots, \lambda_\ell)$. In other words, $\lambda[n]$ is the partition obtained by taking $\lambda$ and adding a large first part to make it have size $n$.  Explicitly, the stability result of Church, Ellenberg, and Farb can be stated as follows:

\begin{theorem}{\textbf{(\cite{CEF} Proposition 3.3.3)}}\label{stab}
Let $V$ be a finitely generated $FI$-module over a field of characteristic zero.  There exist non-negative integer constants $c_\lambda$ independent of $n$ and nonzero for finitely many partitions such that for all $n \gg 0$

$$V_n = \bigoplus_\lambda  c_\lambda S^{\lambda[n]}$$
as a representation of the symmetric group $S_n$, where $S^{\lambda[n]}$ denotes the irreducible Specht module associated to the partition $\lambda[n]$.
\end{theorem}

In positive characteristic, representations of symmetric groups do not in general decompose into a direct sum of irreducible representations. So one might expect that for $FI$-modules over a field $k$ of characteristic $p > 0$ things are more complicated, and indeed that is the case.

 As an example, consider the natural $FI$-module $V$ sending a finite set $S$ to $k[S]$, the space of formal linear combinations of elements of $S$. In this example we have that $V_n \cong k^n$ with the action of $S_n$ permuting the coordinates in the usual way. This is a direct sum of two irreducible representations if $p \nmid n$, and is an indecomposable representation with three irreducible composition factors whenever $p\mid n$ (and $n >2$). So even in this basic example, we see that things are more complicated in positive characteristic.

\end{subsection}

\begin{subsection}{Specht modules and Specht stability}

Specht modules are a well behaved class of symmetric group representations defined over arbitrary commutative rings. While in general they are not irreducible, they play an important role in the representation theory of symmetric groups. We will briefly review some facts about Specht modules that will be important for our purposes. All of the following results can be found in James's book on the representation theory of symmetric groups \cite{James}.


\begin{itemize}

\item For partition $\lambda$ of $n$ the integral Specht module $S^\lambda_\mathbb{Z}$ is a representation of $S_n$ over $\mathbb{Z}$ which is a free and of finite rank as a $\mathbb{Z}$-module.  For an arbitrary commutative ring $k$ the Specht module $S^\lambda_k$ over $k$ is isomorphic to $S^\lambda_\mathbb{Z} \otimes_\mathbb{Z} k$.  From now on we will drop the subscript when the base ring is clear.

\item $S^\lambda_\mathbb{Z}$ is naturally equipped with a nondegenerate $S_n$-equivariant symmetric bilinear form.  This gives rise to a (possibly degenerate) $S_n$-equivariant symmetric bilinear form on $S^\lambda_k$ over any ring $k$.  Let $S^{\lambda \perp}_k$ denote the radical of this form, which is naturally a  $S_n$-invariant subspace of $S^\lambda_k$.

\item Over a field $k$ of characteristic $p$ the quotient $S^\lambda / S^{\lambda \perp}$ is non-zero if and only if $\lambda$ is $p$-regular (meaning $\lambda$ does not have $p$ parts of the same size), in which case it is irreducible.  For $p$-regular partitions $\lambda$, let $D^\lambda$ denote this irreducible quotient. 

\item Over a field of characteristic $p$ the $D^\lambda$ form a complete set of pairwise non-isomorphic irreducible representations of $S_n$ as $\lambda$ runs over the set of all $p$-regular partitions of $n$.

\item Over a field of characteristic $p$ the irreducible composition factors of $S^\lambda$ are all of the form $D^\mu$ with $\mu \ge \lambda$ in the dominance order.  Moreover, for $p$-regular $\lambda$, $D^\lambda$ occurs in $S^\lambda$ with multiplicity one.

\end{itemize}

Specht modules are much better behaved and easier to work with than irreducible representations, so one could hope to generalize representation stability in terms of them.  Putman defined a notion of \emph{Specht stability} \cite{Putman} for a sequence of symmetric group representations in which for all sufficiently large $n$ the $n$th term $V_n$ admits a filtration $$ 0 = V_n^{0} \subset V_n^1 \subset V_n^2 \subset \dots \subset V_n^d = V_n$$ where the graded pieces $V_n^i / V_n^{i-1}$ are isomorphic to Specht modules $S^{\lambda_i[n]}$ with $\lambda_i$ not depending on $n$. 

The previous example sending a set $S$ to $k[S]$ can easily be seen to be Specht stable (of length $2$) by letting $V_n^1$ be the space of formal linear combinations of elements of $S$ where the coefficients sum to zero.  Putman showed that any $FI$-module presented in sufficiently small degree relative to the characteristic of $k$ will exhibit Specht stability (\cite{Putman} Theorem E).  However without the restriction on presentation degree Putman's theorem fails and there are $FI$-modules which do not exhibit Specht stability.

\end{subsection}

\begin{subsection}{Grothendieck groups}

Recall that the Grothendieck group $G_0(\C)$ of an (essentially small) abelian category $\C$ is the abelian group generated by symbols $[X]$ for every object $X$ of $\C$, subject to the relation $[X] - [Y] + [Z] = 0$ for every short exact sequence $$0 \to X \to Y \to Z \to 0$$
of objects in $\C$. If every object of $\C$ has finite length (i.e. if $\C$ is a Krull-Schmidt category) then $G_0(\C)$ is the free abelian group generated by irreducible classes $[X]$, where $X$ is a simple object in $\C$. In this case, given an arbitrary object $Y$, its class $[Y]$ in $G_0(\C)$ is a formal linear combination of irreducible classes $[X]$ with coefficients corresponding to the multiplicity of $X$ in a Jordan-H\"older series for $Y$.

One could ask if for a finitely generated $FI$-module $V$ if the expression for $[V_n]$ in terms of irreducible classes $[D^\lambda]$ in the Grothendieck group of the category of finite dimensional representations of $S_n$ over $k$ stabilizes if we identify $[D^{\lambda[n]}]$ for different values of $n$.  However in our example from before of sending a finite set $S$ to $k[S]$ we have the following inside the Grothendieck group for all $n > 2$:

\[
[V_n]=
\begin{cases}
\hphantom{2}[D^{(n)}] + [D^{(n-1,1 )}]   &\text{if $p \nmid n$},\\[2ex]
 2[D^{(n)}] + [D^{(n-1,1 )}]  &\text{if $p \mid n$}.
\end{cases}
\]

 This example illustrates a general fact about finitely generated $FI$-modules in characteristic $p$. Rather than stabilizing like in characteristic zero, the expression for $[V_n]$ in terms of irreducible classes becomes periodic with period a power of $p$ (see \cite{Harman} section 3.2). While this sort of periodicity is certainly interesting, one often still wants to think of the sequence of symmetric group representations coming from a finitely generated $FI$-module as being stable in some sense. 
 
 \end{subsection}
 
 \begin{subsection}{Virtual Specht stability}
 
The main new idea in this paper is to combine the two ideas above and work in terms of Specht modules inside the Grothendieck groups to obtain a version of representation stability for finitely generated $FI$-modules over fields of arbitrary characteristic. Our main result is the following virtual relaxation of Theorem \ref{stab}, which says that we see stability when expressing $[V_n]$ in terms of Specht classes in the Grothendieck group:

\begin{theorem}\label{virtstab}
Let $V$ be a finitely generated $FI$-module over a field $k$ of positive characteristic.  There exist (possibly negative) integer constants $c_\lambda$ independent of $n$ and nonzero for only finitely many partitions such that for all $n \gg 0$

$$[V_n] = \sum_\lambda  c_\lambda [S^{\lambda[n]}]$$
inside the Grothendieck group of the category of finite dimensional representations of $S_n$ over $k$, where $S^{\lambda[n]}$ denotes the Specht module associated to the partition $\lambda[n]$.
\end{theorem}

 We say that such a sequence of symmetric group representations $V_n$ (whether they come from an $FI$-module or not) satisfying the conclusion of Theorem \ref{virtstab} exhibits \emph{virtual Specht stability}.  In particular, any sequence of representations exhibiting Putman's notion of Specht stability clearly also exhibits virtual Specht stability.

We'll note that in positive characteristics the Specht modules are in general not irreducible, and moreover their classes do not even form a basis for the Grothendieck group (they do span, but there are linear relations between them).  So it is possible that each $[V_n]$ can be expressed in terms of Specht classes in multiple ways.  Virtual Specht stability a priori just requires that among these choices of ways to write $[V_n]$ in terms of Specht classes there is a consistent way to do it for sufficiently large $n$.

Moreover we'll note that for a fixed $FI$-module the choice of the coefficients $c_\lambda$ in the theorem are not unique.  So while this definition of virtual Specht stability will prove easy to work with for our purposes, we'd like to be able to make a statement that involves fewer choices and has a unique output.

Instead of using all Specht modules, we could consider just those corresponding to $p$-regular partitions $\lambda$, which do form a basis of the Grothendieck group (unitriangular relative to the basis of irreducibles).  Instead of asking for the existence of some consistent way of expressing $[V_n]$ in terms of Specht classes, we could ask the stronger question of whether the unique expression for $[V_n]$ as a linear combination of $p$-regular Specht classes stabilizes. Our next theorem says that, in fact, this is equivalent to our (seemingly) weaker notion of virtual Specht stability. Note that $\lambda[n]$ is $p$-regular for $n > 2|\lambda|$ if and only if $\lambda$ is $p$-regular. 

\begin{theorem} \label{unique}

A sequence of representations $V_n$ of $S_n$ over a field of characteristic $p$ exhibits virtual Specht stability if and only if there exist (possibly negative) integer constants $b_\lambda$ independent of $n$ and nonzero for only finitely many $p$-regular partitions $\lambda$ such that for all $n \gg 0$

$$[V_n] = \sum_\lambda  b_\lambda [S^{\lambda[n]}]$$
inside the Grothendieck group of the category of finite dimensional representations of $S_n$ over $k$.

\end{theorem}

So in particular, finitely generated $FI$-modules over a field of characteristic $p$ also exhibit this ``stronger" version of virtual Specht stability.  In fact we will prove a slightly stronger version of Theorem \ref{unique} which also gives control over how two equivalent stable expressions can differ, but we will wait until section \ref{sectequiv} to state the stronger result.

\end{subsection}

\end{section}

\begin{section}*{Acknowledgments}
This result was formulated, proved, and subsequently presented by the author while at the American Institute of Mathematics (AIM) workshop on representation stability in June 2016.  Thanks to Andrew Putman, Steven Sam, Andrew Snowden, David Speyer, and the AIM staff for organizing the workshop and for inviting me to speak. Thanks as well to Pavel Etingof, Benson Farb, Andrew Putman, and Jordan Ellenberg for helpful comments and conversations. This work was partially supported by the National Science Foundation Graduate Research Fellowship under Grant No. 1122374.

\end{section}

\begin{section}{Virtual Specht stability for $FI$-modules}

In this section we will prove Theorem \ref{virtstab}. The main tool we will need is a powerful structural result for $FI$-modules over arbitrary Noetherian rings due to Nagpal. Before stating the result, we will review a few definitions.

Let $W$ be a representation of $S_m$ over an arbitrary noetherian ring $k$ (for our purposes, $k$ will be a field of positive characteristic). The $FI$-module $M(W)$ is the functor that takes a finite set $S$ to the vector space $k[\text{Hom}_{FI}([m],S)] \otimes_{k[S_m]} W$.  In other words, $M(W)_n$ is the representation of $S_n$ obtained by first extending $W$ to a $S_m \times S_{n-m}$ representation by letting the second factor act trivially, and then inducing up to $S_n$.

Such modules, and direct sums thereof, are exactly those $FI$-modules which admit the structure of an $FI\sharp$-module (See \cite{CEF}, section 4).  We say that an $FI$-module $V$ is $\sharp$-filtered if it admits a filtration $$0 = V^0 \subset V^1 \subset \ldots \subset V^d = V$$ of $FI$-submodules such that the graded pieces $V^i / V^{i-1}$ are of the form $M(W)$ for some $W$.  Informally, the result of Nagpal we will need says that, up to torsion, every $FI$-module admits a resolution by $\sharp$-filtered modules. More precisely, it says:

\begin{theorem}{\textbf{(\cite{Nagpal} Theorem A)}}\label{rohit} For an arbitrary finitely generated $FI$-module $V$, there exist $\sharp$-filtered $FI$-modules $J^1, J^2, J^3, \dots, J^N$ and maps of $FI$-modules $\phi^0: V \to J^1, \ \phi^1: J^1 \to J^2, \ \dots, \ \phi^{N-1}: J^{N-1} \to J^N$ such that for all $n \gg 0$

$$0 \to V_n \to J^1_n \to J^2_n \to \dots \to J^N_n \to 0$$
is an exact sequence of $S_n$ representations.

\end{theorem}

In addition to Nagpal's result on $FI$-modules, we will need one standard fact from the modular representation theory of symmetric groups. Informally, the following theorem says that at the level of Grothendieck groups induction between Specht modules behaves the same as in characteristic zero.  

\begin{theorem} \textbf{(\cite{James} Corollary 17.14)} \label{spechtfilt} For arbitrary partitions $\lambda \vdash n, \ \mu \vdash m$ and over an arbitrary field $k$, the induced representation $V =  Ind_{S_n \times S_m}^{S_{n+m}}(S^\lambda \boxtimes S^\mu)$ has a filtration $$0 = V^0 \subset V^1 \subset \ldots \subset V^d = V$$ such that the graded pieces $V^i / V^{i-1}$ are isomorphic to Specht modules $S^\nu$, with $S^\nu$ occurring with multiplicity equal to the Littlewood-Richardson coefficient $c_{\lambda,\mu}^{\nu}$.

\end{theorem}

\medskip

\noindent \textbf{Proof of Theorem \ref{virtstab}:}  The proof will be by a series of reductions to a known case from characteristic zero. First, Nagpal's result allows us to replace $[V_n]$ by an alternating sum of the $[J^i_n]$'s inside the Grothendieck group for all sufficiently large $n$. Hence it is enough to prove the result for $\sharp$-filtered $FI$-modules.

Next, we know by definition that $\sharp$-filtered modules admit a filtration with quotients of the form $M(W)$, where $W$ is a representation of some $S_m$. At the level of Grothendieck groups such filtrations just become sums, so it is enough to prove the result for the modules of the form $M(W)$.  

Now since the functor $W \to M(W)_n$ (which again is given by extending $W$ to a $S_m \times S_{n-m}$ representation by letting $S_{n-m}$ act trivially and then inducing up to $S_n$) is exact, it descends to a linear map on Grothendieck groups, so it is enough to prove the claim for a collection of representations $W$ that span the Grothendieck groups of symmetric groups.

We know that the Specht modules $S^\lambda$ are exactly such a class of representations
.  Finally, Theorem \ref{spechtfilt} tells us that $M(S^\lambda)_n = Ind_{S_m \times S_{n-m}}^{S_{n}}(S^\lambda \boxtimes S^{(n-m)})$ has a filtration by Specht modules with multiplicities the same as in characteristic zero, where we know stabilization occurs. $\square$

\medskip
One consequence of representation stability in characteristic zero is that the characters $\chi_{V_n}$ of the $S_n$ representations $V_n$ are eventually polynomial functions in the number of cycles of given lengths (\cite{CEF} Theorem 1.5). Since Specht modules are defined over the integers their Brauer characters agree with their usual characters in characteristic zero, this along with Theorem \ref{virtstab} immediately implies the following positive characteristic analog of this fact:

\begin{corollary}
If $V$ is a finitely generated $FI$-module over a field of positive characteristic then the sequence of Brauer characters $\hat{\chi}_{V_n}$ of the $S_n$ representations $V_n$ is eventually polynomial.
\end{corollary}

\end{section}
 
\begin{section}{Equivalent stable presentations}\label{sectequiv}

For this section all representations are over a fixed field $k$ of characteristic $p > 0$. We will refer to ``the Grothendieck group for $S_n$" which should be taken to mean ``the Grothendieck group of the category of finite dimensional representations of $S_n$ over $k$".
 
For a fixed value of $n$ we know that the Specht classes $[S^\lambda]$ span the Grothendieck group for $S_n$, but since there more Specht modules than irreducible representations there are linear relations among the Specht classes.  
 
One might hope that when we pass to the stable setting that these relations go away, since after all we are requiring a relation to hold in infinitely many Grothendieck groups simultaneously.  Unfortunately, there is an easy source of stable linear relations between Specht classes that we will outline now.

Let  $\sum c_\lambda[S^\lambda] = 0$ be a linear relation between Specht classes inside the Grothendieck group for $S_m$.  Then we know that $$\sum c_\lambda [M(S^\lambda)_n] = \sum c_\lambda [Ind_{S_m \times S_{n-m}}^{S_{n}}(S^\lambda \boxtimes S^{(n-m)})] =0$$
 in the Grothendieck group for $S_n$ for all $n > m$.  If we then expand each term into Specht classes using Theorem \ref{spechtfilt} we then obtain a stable expression of the form $\sum d_\lambda[S^{\lambda[n]}] = 0$ with $d_\lambda$ not depending on $n$, $d_\lambda = c_\lambda$ if $|\lambda| = m$, and $d_\lambda =0$ if $|\lambda|>m$ which holds for all $n > 2m$.  We will call such an expression an $\emph{induced expression for zero}$.
 
In particular this implies that the stable expression for a finitely generated $FI$-module in terms of Specht classes guaranteed by Theorem \ref{virtstab} is far from unique, we can just add an induced expression for zero as constructed above to obtain another equivalent stable expression.
 
We'd like to uniqueness statement about the expressions in terms of Specht classes, and there are (at least) two ways we could try to do this. First, we can restrict to those Specht classes corresponding to $p$-regular $\lambda$ as in the statement of Theorem \ref{unique}, and ask if the unique expression in terms of these stabilizes.  Alternatively, we could ask if a stable expression for an $FI$-module (or other virtually Specht stable sequence) is unique up to a linear combination of induced expressions for zero.  

The following proposition does both of these simultaneously, immediately implying both Theorem \ref{unique} and that any two stable expressions for the same $FI$-module differ by a linear combination of induced expressions for zero.

\begin{proposition}\label{unique2}
For any expression of the form $\sum c_\lambda[S^{\lambda[n]}]$ there is an expression $\sum d_\lambda[S^{\lambda[n]}]$ equivalent to the first inside the Grothendieck group for $S_n$ for all sufficiently large $n$ such that $d_\lambda = 0$ for $p$-singular $\lambda$, and the difference $\sum c_\lambda[S^{\lambda[n]}] - \sum d_\lambda[S^{\lambda[n]}]$ is a linear combination of induced expressions for zero.

\end{proposition}

\noindent \textbf{Proof:}  Let $\lambda_0$ be such that $\lambda_0[n]$ is minimal in the dominance order among those $p$-singular partitions $\lambda$ with $c_\lambda \ne 0$ and  let $m = |\lambda_0|$. iInside the Grothendieck group for $S_{m}$ we know we can write:
$$[S^{\lambda_0}] = \sum_{\lambda > \lambda_0 } b_\lambda [S^\lambda]$$ 
for some integer coefficients $b^\lambda$.  Applying the functor $M(*)_n$ to both sides, expanding in terms of Specht classes using Theorem \ref{spechtfilt}, and rearranging gives an induced expression for zero of the form $$[S^{\lambda_0[n]}] - \sum_{\lambda[n] > \lambda_0[n] } b'_\lambda [S^{\lambda[n]}] = 0$$
for all sufficiently large $n$.  If we subtract off $c_{\lambda_0}$ times this from our expression $\sum c_\lambda[S^{\lambda[n]}]$  we get an equivalent expression $\sum c'_\lambda[S^{\lambda[n]}]$ where $c'_{\lambda_0} = 0$ differing from the first by Specht classes corresponding to partitions $\lambda[n]$ larger than $\lambda_0[n]$ in the dominance order. 

Since the (stable) dominance order has no infinite ascending chains, we can repeat this process until there are no non-zero coefficients for Specht classes corresponding to $p$-singular partitions. $\square$

\end{section}
 
\begin{section}{Examples}
Finally we will finish by giving a couple of examples to illustrate our results.

\medskip

\noindent \textbf{Example 1:} Returning to our example from the first section of the $FI$-module $V$ that sends a finite set $S$ to $k[S]$.  Inside $V_n$ there is the subspace of formal linear combinations of elements of $[n]$ where the sum of the coefficients is zero. This space is isomorphic to the Specht module $S^{(n-1,1)}$.  The quotient of $V_n$ by this subspace is a copy of the trivial representation, which itself is the Specht module $S^{(n)}$.  Hence we see that $[V_n] = [S^{(n-1,1)}]+[S^{(n)}]$ for all $n \ge 2$.   This holds even when the characteristic of $k$ divides $n$, in which case the Specht module $S^{(n-1,1)}$ is not irreducible as it contains a $1$-dimensional invariant subspace. Of course this example also satisfies Putman's stronger notion of (non-virtual) Specht stability.

\medskip

\noindent \textbf{Example 2:} Let $k$ be a field of characteristic $5$, and consider the standard representation of $S_5$ acting on $k^5$ by permuting the coordinates.  This contains two proper $S_5$ invariant subspaces: the $4$-dimensional space of vectors with coordinate sum zero (i.e. the Specht module $S^{(4,1)}$), and the $1$-dimensional space of invariants spanned by the vector $(1,1,1,1,1)$.   Since we are in characteristic $5$, this vector $(1,1,1,1,1)$ has coordinate sum zero and therefore the space of invariants is contained in the Specht module.  Let $W =  D^{(4,1)}$ be the $3$-dimensional quotient of $S^{(4,1)}$ by this invariant subspace.

Now consider the $FI$-module $M(W)$.  Unlike the first example, the representations $M(W)_n$ in general do not admit a filtration where the successive quotients are Specht modules.  Nevertheless, by construction we see that in the Grothendieck group $[W] = [S^{(4,1)}] - [S^{(5)}]$.  Applying the Pieri rule to this expression term wise and simplifying we obtain an expression for $[M(W)_n]$ for all $n \ge 10$: $$[M(W)_n] = [S^{(n-5,4,1)}] + [S^{(n-4,3,1)}] + [S^{(n-3,2,1)}] + [S^{(n-2,1,1)}] - [S^{(n-5,5)}]$$

\noindent \textbf{Example 3:} Let $W = S^{(1,1,1)}$ be the sign representation of $S_3$ over a field of characteristic $3$.  We know that $M(W)$ exhibits Specht stability by Theorem \ref{spechtfilt}, and at the level of Grothendieck groups we have that $$[M(W)_n]= [S^{(n-2,1,1)}]+[S^{(n-3,1,1,1)}]$$
for $n \ge 4$. This satisfies our definition of virtual Specht stability, however it involves a $3$-singular partition. Theorem \ref{unique} tells us that we should also see stability when expressing $[M(W)_n]$ in terms of Specht classes for $3$-regular partitions.  In this case we have that $$[S^{(1,1,1)}] = [D^{(2,1)}] = [S^{(2,1)}]-[S^{(3)}]$$ which by the Pieri rule tells us that
$$[M(W)_n] = [S^{(n-3,2,1)}]+[S^{(n-2,1,1)}] - [S^{(n-3,3)}]$$

\end{section}

\end{document}